\documentclass[preprint,11pt]{elsarticle}

\usepackage[pagewise]{lineno}\linenumbers

\usepackage[utf8]{inputenc}

\pdfoutput=1

\usepackage{mdwlist}
\usepackage{graphicx}
\usepackage{amsfonts}
\usepackage{ulem}
\usepackage{amsmath}
\usepackage{amsthm}
\usepackage{amssymb}
\usepackage{fancyhdr}
\usepackage{titlesec}
\usepackage{indentfirst}
\usepackage{booktabs}
\usepackage{verbatim}
\usepackage{color}
\usepackage{latexsym}
\usepackage{amscd}
\usepackage{esvect}
\usepackage{graphicx}
\usepackage{upgreek}
\usepackage{subcaption}   
\usepackage{caption}
\usepackage[page,header]{appendix}
\usepackage{titletoc}
\usepackage{mathrsfs}
\usepackage[numbers]{natbib}
\usepackage{float}
\usepackage{lipsum}
\usepackage[export]{adjustbox} 
\numberwithin{equation}{section}
\newtheorem{theorem}{Theorem}[section]
\newtheorem{lemma}[theorem]{Lemma}

\newtheorem{remark}[theorem]{Remark}

\allowdisplaybreaks
\usepackage{todonotes}
\usepackage[colorlinks,citecolor=red,urlcolor=blue,bookmarks=false,hypertexnames=true]{hyperref} 

\newcommand{\beq}{\begin{equation}}
\newcommand{\eeq}{\end{equation}}
\newcommand{\ben}{\begin{eqnarray}}
\newcommand{\een}{\end{eqnarray}}
\newcommand{\beno}{\begin{eqnarray*}}
\newcommand{\eeno}{\end{eqnarray*}}

\journal{Journal of Discrete and Continuous Dynamical System – S}
\begin{document}

\begin{frontmatter}

\title{Solving the linear transport equation by a deep neural network approach}

\author[label1]{Zheng Chen}
\address[label1]{Department of Mathematics, University of Massachusetts, Dartmouth, MA, 02747, zchen2@umassd.edu}
\author[label2]{Liu Liu}
\address[label2]{Department of Mathematics, The Chinese University of Hong Kong, Hong Kong, lliu@math.cuhk.edu.hk}
\author[label3]{Lin Mu}
\address[label3]{Department of Mathematics, University of Georgia, Athens, GA 30602, linmu@uga.edu}

\begin{abstract}
In this paper, we study the linear transport model by adopting the {\it deep learning method}, in particular the deep neural network (DNN) approach. While the interest of using DNN to study partial differential equations is arising,  
here we adapt it to study kinetic models, in particular the linear transport model. Moreover, theoretical analysis on the convergence of the neural network and its approximated solution towards the analytic solution is shown. We demonstrate the accuracy and effectiveness of the proposed DNN method in the numerical experiments. 
\end{abstract}


\end{frontmatter}



\section{Introduction}

The deep learning method is a new approach for solving partial differential equations (PDEs) that can resolve some difficulties appeared in traditional numerical methods \cite{Karn-NN}, 
such as the expensive computational cost for high-dimensional problems, the challenges dealing with complex boundary conditions, truncation of the velocity domain (in some models such as the Boltzmann equation \cite{Cercignani-Book}, the velocity lies in the whole space thus a truncation in numerical discretization is necessary and one needs to study decaying properties of the solution). One of the deep learning approaches that worth mentioning is the physical informed neural networks (PINN) \cite{PINN-1,PINN-2,PINN-UQ}. This method incorporates the easy-to-use auto-differentiation in the current software and the physics information of the PDEs under study into the network, by constructing total loss functions that involves the PDE residuals and other constraints such as the initial and boundary conditions. 


Besides, the deep learning algorithm, as a mesh-free method, has the advantage of being intuitive and easy to be executed. For example, instead of designing mass (or momentum and energy) conservative schemes for kinetic models, which can be challenging for traditional numerical methods and may need extra efforts \cite{MM-Cons}, one can simply involve the derivative in time of conserved quantities of interests in the total loss function, as was done in \cite{Hwang} for the kinetic Fokker-Planck equation. 
However, we mention that there exist indeed some weaknesses of the deep learning approach.
First, there is no guarantee that the deep learning algorithms will converge and it is theoretically difficult to show their convergence in practice. It is also hard to evaluate the accuracy of the deep neural network (DNN) approach in contrast with conventional numerical methods such as finite volume or finite element methods.  

We mention some other advantages of using the DNN approach to solve our kinetic problems, in particular the linear transport model: 
a) to obtain the distribution function at any given time $t$, position $x$, velocity $\Omega$, instead of only the discrete values at uniform mesh in traditional finite difference numerical methods; as a mesh-free method, works efficiently for high-dimensional physical space problems; 
b) to avoid high computational cost on the simulation of kinetic equations due to the spacial and velocity variables, and the integral-based, nonlocal collision operators. 

Many kinetic applications are modeled by a linear kinetic transport equation that describes how kinetic particles get collisions and absorption through a material medium and evolve in time. This model has been applied in a wide variety of fields, such as atmosphere and ocean modeling \cite{coakley2014atmospheric,stamnes2017radiative,zdunkowski2007radiation}, astrophysics \cite{peraiah2002introduction}, and neutron transport and nuclear physics \cite{case1967linear}. 
Understanding such models both theoretically and numerically are very important. 
Many numerical methods have been developed to simulate model equations, such as Monte Carlo methods, discrete ordinate methods \cite{alcouffe1985first,hauck2019filtered,lathrop1968ray}, and moment methods \cite{chen2019multiscale,garrett2013comparison,hauck2010positive,laiu2016positive,mcclarren2010robust}. The collisions induce some mathematical structure, which is utilized to design many algorithms.   
There are domain decomposition methods \cite{bourgat1996coupling,golse2003domain,klar1994domain,klar2000transition}, perturbative methods \cite{brunner1999collisional,degond2005smoothdiffusion,degond2005smoothhydrodynamic,degond2006macroscopic,dimarco2008hybrid,dimarco2010fluid,liu2004boltzmann,parker1993fully}, asymptotic preserving (AP) numerical methods \cite{chen2017dg,jang2015high,jin1999efficient,jin2010asymptotic}, and collision-based hybrid algorithms \cite{crockatt2017arbitrary,hauck2013collision,heningburg2020hybrid}. 
The equation describes particle advection and particle interactions, which are basic features shared by other kinetic models such as the ones describe dilute gases  \cite{Cercignani-1988, Cercignani-Illner-Pulvirenti-1994,chapman1970mathematical}; neutron \cite{ case1967linear, dautray2012mathematical,davison1957neutron,lewis1984computational}, photon \cite{mihalas1999foundations,Pomraning-1973}, and neutrino \cite{mezzacappa1999neutrino} radiation; charged transport in semiconductor devices \cite{ decaria2020analysis,Markovich-Ringhofer-Schmeiser-1998, selberherr2012analysis}; and ionized plasma \cite{boyd2003physics,hazeltine2004framework}.   
Successful development of DNN for linear kinetic equation will indicate good potential of DNN approach to simulate more models listed above.

The trend of using DNN to solve PDE problems is arising, we refer to \cite{Hwang} for a review of the literature and omit it here. In this work we have made our main contribution by adapting it to study the kinetic problems, in particular the linear transport model with practical applications. 
In the theoretical proof, we show that as the number of neurons goes to infinity, 
1) the total loss function goes to zero; 
2) the neural network solution converges point-wisely to a {\it priori} analytic solution of the linear transport equation; 
We use energy estimates analysis to study the convergence of the neural network solution to the a priori analytic solution in this work.

This paper is organized by the following. In Section \ref{sec:LTE}, we introduce the linear transport equation under study. 
In Section \ref{sec:NN}, we review and discuss the neural network framework and method for solving general PDEs. 
Then, we give main convergence results in Section \ref{sec:Anal}, which show that 
1) the loss function {goes to zero as the neural network converges};  
2) the neural network solution converges point-wisely to the analytic solution when the loss function converges to zero. 
In Section \ref{sec:Num}, numerical simulation of the DNN approach of a few examples with practical applications for the linear transport model will be presented, and we shall demonstrate the efficiency and accuracy of the proposed DNN method. 
Finally, we make conclusions and mention some future work in Section \ref{sec:Con}. 

\section{The transport equation}
\label{sec:LTE}
In this paper, we consider a single-group, linear transport equation of the form
\begin{equation}\label{eqn:3Dtransport}
    \frac{1}{c}\frac{\partial \psi}{\partial t} + \Omega \cdot \nabla_x \psi + \sigma_t \psi = \sigma_s \mathcal{S} \psi + \mathcal{Q},
\end{equation}
where $x \in \mathrm{X} \subset \mathbb{R}^3$ is position, $\Omega \in \mathbb{S}^2$ is direction of flight, $t>0$ is time.  The particle speed $c$ is a fixed scalar and $\mathcal{Q} = \mathcal{Q}(x)$ is a source.  The material is characterised by the absorption, scattering, and total cross-sections, denoted as $\sigma_a$, $\sigma_s$ and $\sigma_t = \sigma_a + \sigma_s$, respectively.

For each fixed $x$, the scattering operator $\mathcal{S}: L^2(\mathbb{S}^2) \rightarrow L^2(\mathbb{S}^2)$ is a bounded linear operator, and for any $f(x,\cdot) \in L^2(\mathbb{S}^2)$,
\begin{equation}
    \left(\mathcal{S} f\right)(x,\Omega) = \int_{\mathbb{S}^2}\, g(x,\Omega\cdot\Omega') f(x,\Omega') \, d\Omega',
\end{equation}
where $g$ is the scattering angular redistribution function and is normalized to one, i.e. for any $x \in \mathrm{X}$,
\begin{equation}
    \int_{\mathbb{S}^2}\, g(x,\Omega\cdot\Omega') \, d\Omega' = 2\pi \int_{-1}^1 \, g(x,\mu)\, d\mu = 1.
\end{equation}
\section{The Neural Network approach}
\label{sec:NN}


We review the deep neural network (DNN) structure and approach for general PDEs. Please see Figure \ref{Fig:NN} below. 
Denote the approximated solution to the linear transport model \eqref{eqn:3Dtransport} by 
$\psi^{nn}(t,x,\Omega;m,w,b)$ and suppose the neutral network has $L$ layers; the input layer takes $(t,x,\Omega)$ and the final layer gives $\psi^{nn}(t,x,\Omega;m,w,b)$ as the output. The relation between the $l$-th and $(l+1)$-th layer ($l=1,2, \cdots L-1)$ is 
given by
\begin{equation}\label{two-layer}
    \theta_j^{(l+1)}=\sum_{i=1}^{m_{l}} w_{j i}^{(l+1)} \bar{\sigma}_{l} (\theta_{i}^{l})+b_{j}^{(l+1)},
\end{equation} 
where $m=\left(m_{0}, m_{1}, m_{2}, \dots, m_{L-1}\right)$, $w=\left\{w_{j i}^{(k)}\right\}_{i, j, k=1}^{m_{k-1}, m_{k}, L}$ and $b=\left\{b_{j}^{(k)}\right\}_{j=1, k=1}^{m_{k}, L}$. More specifically, 
\begin{itemize}
    \item $\theta_i^l$: the $i$-th neuron in the $l$-th layer 
    \item $\bar{\sigma}_l$: the activation function in the $l$-th layer
    \item $w_{j i}^{(l+1)}$: the weight between the $i$-th neuron in the $l$-th layer and the $j$-th neuron in the $(l+1)$-th layer
    \item $b_j^{(l+1)}$: the bias of the $j$-th neuron in the $(l+1)$-th layer
    \item $m_l$: the number of neurons in the $l$-th layer. 
\end{itemize}

\begin{figure}[!h]
\centering
\includegraphics[width=.6\textwidth,height = .4\textwidth]{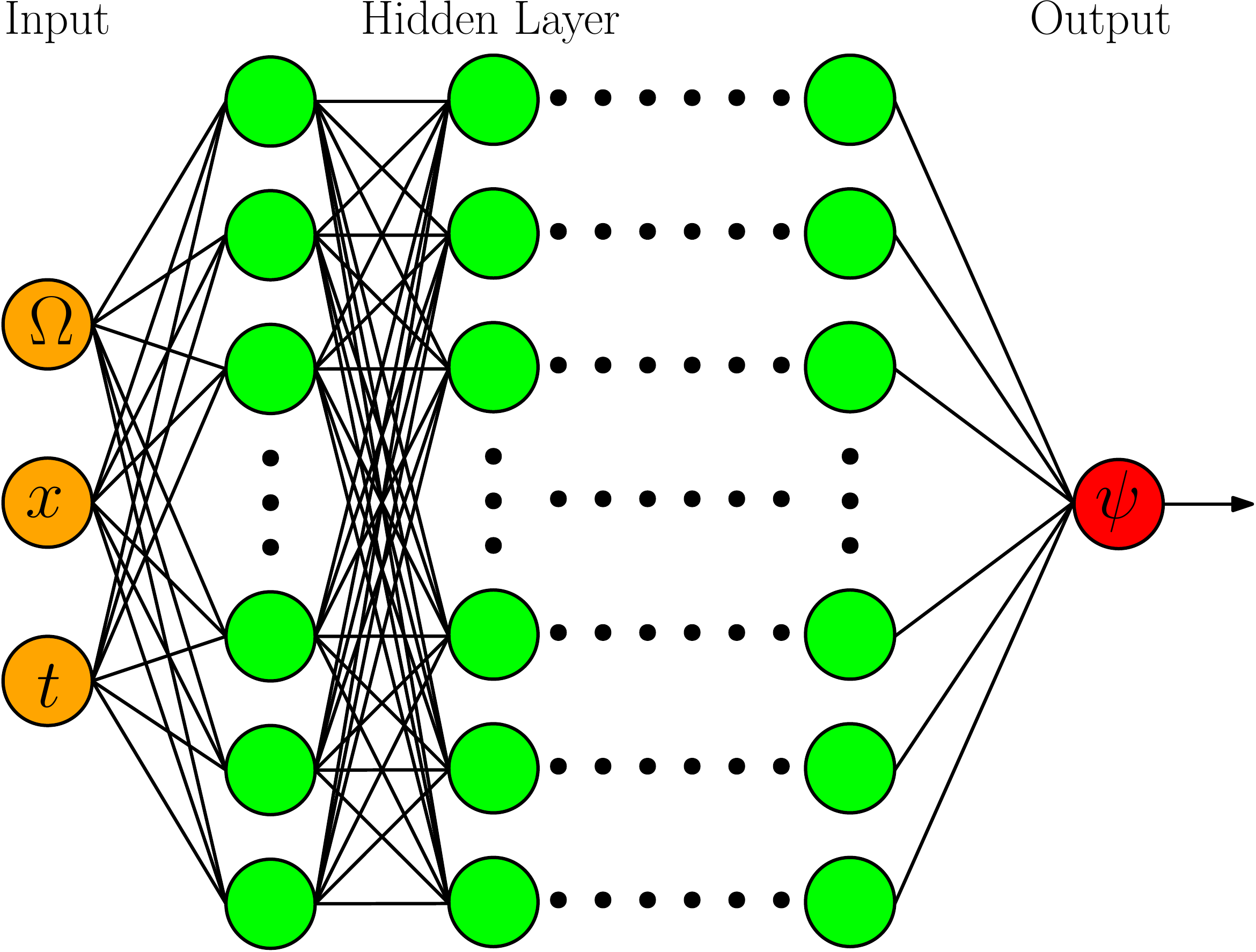}
\caption{Structure of the deep neural network.}
\label{Fig:NN}
\end{figure}
\subsection{Definition of loss functions}

Recall the linear transport equation \eqref{eqn:3Dtransport}, for simplicity, we let $c=1$, $\mathcal Q=0$ and denote the operator
$\mathcal {L}(\psi) = \sigma_s \mathcal{S}\psi - \sigma_t \psi$.
By \eqref{eqn:3Dtransport}, we consider the linear transport equation given by
\begin{eqnarray}
\label{LTE}
\begin{cases}
&\frac{\partial \psi}{\partial t} + \Omega \cdot \nabla_x \psi = \mathcal{L}(\psi), 
 \\[4pt]
& \psi(t=0,x,\Omega) = \psi_0(x,\Omega).
\end{cases}
\end{eqnarray}
The loss function for the governing equation \eqref{LTE} is defined by
\begin{align}
\label{Loss1}
\begin{split}
\displaystyle Loss_{\text{GE}} &\displaystyle  := \int_0^T \int_{\mathrm{X}}\int_{\mathbb S^2} \Big| \partial_t \psi^{nn}(t,x,\Omega;m,w,b) + \Omega \cdot \nabla_x \psi^{nn}(t,x,\Omega;m,w,b) \\[4pt]
&\displaystyle\quad - \mathcal L(\psi^{nn})(t,x,\Omega;m,w,b) \Big|^2\,  d\Omega\, dx\, dt \\[4pt]
&\displaystyle \approx \sum_{i,j,k} \Big| \partial_t \psi^{nn}(t_i, x_j, \Omega_k; m,w,b) + \Omega \cdot \nabla_x \psi^{nn}(t_i, x_j, \Omega_k; m,w,b) \\[4pt]
&\displaystyle \quad - \mathcal L(\psi^{nn})(t_i, x_j, \Omega_k; m,w,b)\Big|^2 \omega_{i,j,k}, 
\end{split}
\end{align}
where $\{(t_i, x_j, \Omega_k)\}$ are the discretization points used to approximate the integral over the whole domain, with $\omega_{i,j,k}$ the associated weight for point $(t_i, x_j, \Omega_k)$. For uniform points, it is simply given by 
\begin{equation}
    \omega_{i,j,k} = \frac{1}{N_i N_j N_k},
\end{equation}
where $N_i, N_j, N_k$ are the numbers of discretization points used in $(t,x,\Omega)$ direction respectively.
%
We shall introduce the discretized form 
$\mathcal L(\psi^{nn})(t,x,\Omega; m,w,b)$ in Section \ref{sec:Num} later. 

The loss term for the initial condition is defined by
\begin{align}
\label{Loss2}
\begin{split}
\displaystyle Loss_{\text{IC}} &:= \int_{\mathrm{X}}\int_{\mathbb S^2} 
\Big|\psi^{nn}(0,x,\Omega) - \psi_0(x,\Omega)\Big|^2\, d\Omega\, dx \\[4pt]
\displaystyle & \approx \frac{1}{N_{\text{IC}}}\sum_{(x,\Omega) \in \mathbb{X}_{\text{IC}}}
\Big|\psi^{nn}(0,x,\Omega) - \psi_0(x,\Omega)\Big|^2,
\end{split}
\end{align}
where $\mathbb{X}_{\text{IC}}$ denotes the set of sampling points in $\mathrm{X}\times\mathbb S^2$ used to approximate the initial data with total number $N_{\text{IC}}$. 

Define $n_x$ the {\it unit} outward normal vector on the boundary $\partial\mathrm{X}$, let $\gamma \stackrel{\text { def }} = \partial\mathrm{X}\times \mathbb{S}^2$. 
This phase boundary can be split into an outgoing boundary $\gamma_{+}$, incoming boundary $\gamma_{-}$ and a singular boundary $\gamma_0$, defined by

\begin{align}
\gamma_{+} :=\left\{(x, \Omega) \in \partial\mathrm{X}\times \mathbb{S}^2: n_x \cdot \Omega>0\right\}, \\[2pt]
\gamma_{-} :=\left\{(x, \Omega) \in \partial\mathrm{X}\times \mathbb{S}^2: n_x \cdot\Omega<0\right\}, \\[2pt]
\gamma_{0} :=\left\{(x, \Omega) \in \partial\mathrm{X}\times \mathbb{S}^2: n_x \cdot\Omega=0\right\}. 
\end{align}
%
We assume the inflow boundary condition 
\begin{equation}\label{bc}
    \psi(t,x,\Omega)|_{\gamma^{-}}=h(t,x,\Omega)\quad \text{ for } x\in\partial\mathrm{X},
\end{equation}
with the loss term defined by 
\begin{align}
\label{Loss3}
\begin{split}
\displaystyle Loss_{\text{BC}} &:= \int_0^T \int_{\gamma_{-}} \Big|\psi^{nn}(t,x,\Omega;m,w,b) - h(t,x,\Omega;m,w,b)\Big|^2\, ds\, dt  \\[4pt]
\displaystyle &\approx \frac{1}{N_{\text{BC}}} \sum_{(t,x,\Omega)\in \mathbb{X}_{\text{BC}}} 
\Big|\psi^{nn}(t,x,\Omega;m,w,b) - h(t, x, \Omega;m,w,b)\Big|^2, 
\end{split}
\end{align}
where $\mathbb{X}_{\text{BC}}$ denotes the set of sampling points in $[0,T]\times\gamma_{-}$ used to approximate the inflow boundary data, with total number $N_{\text{BC}}$.
Similarly, for the specular reflection boundary condition given by 
\begin{equation}\label{rbc}
    \psi(t,x,\Omega) = \psi(t,x, \mathcal{R}\Omega) \quad \text{ for  } (t,x,\Omega)\in \gamma, 
\end{equation}
where $\mathcal{R}\Omega = \Omega - 2 (n_x\cdot\Omega)n_x$, one can define the corresponding loss term 
\begin{align}
\label{Loss_Spec}
\begin{split}
\displaystyle \widetilde{Loss}_{\text{BC}} &:= \int_0^T \int_{\gamma} \Big|\psi^{nn}(t,x,\Omega;m,w,b) - \psi^{nn}(t,x,\mathcal{R}\Omega;m,w,b)\Big|^2\, ds\, dt  \\[4pt]
\displaystyle &\approx \frac{1}{\tilde{N}_{\text{BC}}} \sum_{(t,x,\Omega)\in \tilde{\mathbb{X}}_{\text{BC}}} 
\Big|\psi^{nn}(t,x,\Omega;m,w,b) - \psi^{nn}(t,x,\mathcal{R}\Omega;m,w,b)\Big|^2, 
\end{split}
\end{align}
where $\tilde{\mathbb{X}}_{\text{BC}}$ denotes the set of sampling points in $[0,T]\times\gamma$, with total number $\tilde{N}_{\text{BC}}$.

Adding up \eqref{Loss1}, \eqref{Loss2} and \eqref{Loss3} defines the total loss function as follows
\begin{equation} \label{Total-Loss}
{Loss}_{\text{Total}} := Loss_{\text{GE}} + Loss_{\text{IC}} + Loss_{\text{BC}}. 
\end{equation}

\section{Analysis main results}
\label{sec:Anal}


In this section, we show two main theoretical results. We first prove that there exists a sequence of parameters (weights, biases and number of nodes) such that the total loss function, defined earlier in \eqref{Total-Loss}, converges to $0$. 
Consequentially, we then show that the sequence of these neural networks equipped with such parameters converge to the analytic solution of the linear transport model \eqref{eqn:3Dtransport}. 
We first review some preliminary results about the existence of the approximated neural network solution in section \ref{subsec:Pre}. 

\subsection{Preliminaries}
\label{subsec:Pre}

The neural network architecture was first introduced in \cite{NN1943}. Later in \cite{UAT}, Cybenko established sufficient conditions where a continuous function can be approximated by finite linear combinations of single hidden layer neural networks, followed by 
the work in \cite{Multilayer} that extends the theory to the multi-layer network case. 

In the following, we will review an important theorem on the existence of the approximated neural network solution, that will prepare us for the analysis shown later in this section. 
First, we recall the Universal Approximation Theorem (UAT) in \cite{UAT} which uses the two-level neural network. We adapt to discuss our model equation and use the notations introduced in Section \ref{sec:NN}. 
\begin{lemma}
Suppose the solution to \eqref{eqn:3Dtransport} satisfies $\psi \in C^1([0,T]) \cap C^1(\mathrm{X})$. Also, let the activation function $\bar\sigma$ be any non-polynomial function in $C^1(\mathbb R)$. Then for any $\epsilon>0$, there exists a two-layer neural network
$$ \psi^{nn}(t,x,\Omega) = \sum_{i=1}^{m_1} w_{1i}^{(2)} \bar\sigma\left( \left(w_{i 1}^{(1)}, w_{i 2}^{(1)}, w_{i 3}^{(1)}\right) \cdot(t, x, \Omega)+b_{i}^{(1)}\right) + b_{1}^{(2)}, 
$$
such that 
\begin{equation}\label{UAT}
    \left\| \psi - \psi^{nn} \right\|_{L^{\infty}(K)} < \epsilon, \quad 
\left\| \partial_t ( \psi - \psi^{nn}) \right\|_{L^{\infty}(K)} < \epsilon, \quad
\left\| \nabla_x ( \psi - \psi^{nn}) \right\|_{L^{\infty}(K)} < \epsilon, 
\end{equation}
where the domain $K$ denotes $ [0,T]\times\mathrm{X}\times\mathbb{S}^2 $. 
\end{lemma}
\begin{remark}
Notice that the above result can be generalized to the neural network with several hidden layers \cite{Multilayer}. For simplicity, we review the above UAT in the two-layer neural network case. 
\end{remark}

\subsection{Convergence of the loss function}

We first show that a sequence of neural network solutions to \eqref{LTE}
exists such that the total loss function converges to zero, if the analytic solution $\psi \in C^1([0,T]) \cap C^1(\mathrm{X})$.

\begin{theorem}
\label{Thm1}
Assume the solution $\psi$ to \eqref{eqn:3Dtransport} is sufficiently smooth in its physical domain $K=[0,T]\times\mathrm{X}\times\mathbb{S}^2$, in the sense that $\psi \in C^1([0,T]) \cap C^1(\mathrm{X})$. Then there exists a sequence of neural network parameters $\{m_{[j]}, w_{[j]}, b_{[j]}\}_{j=1}^{\infty}$ such that the sequence of DNN solutions with $m_{[j]}$ nodes, denoted by $\{\psi_j(t,x,\Omega) = \psi^{nn}(t,x,\Omega;m_{[j]}, w_{[j]}, b_{[j]})\}_{j=1}^{\infty}$ satisfies 
$$ Loss_{\text{Total}}(\psi_j)\to 0, \qquad \text{as  } j\to\infty. $$
\end{theorem}
\bigskip
{\bf Proof. }  
Define a sequence of small numbers $\epsilon_j = \frac{1}{j}$. 
By the UAT above, for any $\epsilon_j$, there exists a DNN solution
\begin{equation}
\psi_j(t, x, \Omega)=\sum_{i=1}^{m_{[j], 1}} w_{[j], 1 i}^{(2)} \bar{\sigma}\left(\left(w_{[j], i 1}^{(1)}, w_{[j], i 2}^{(1)}, w_{[j], i 3}^{(1)}\right) \cdot(t, x, \Omega)+b_{[j], i}^{(1)}\right)+b_{[j], 1}^{(2)}
\end{equation}
such that 
 $||\psi_j - \psi||_{L^{\infty}(K)} < \epsilon_j$ and similarly for their first-order derivatives in $t$ and $x$ as in \eqref{UAT}. Here the neural network parameters $m$, $w$ and $b$ satisfies the definition in \eqref{two-layer}, except that there is an extra subscript $[j]$, which describes the index of the sequence 
 of DNN solutions $\{\psi_j \}_{j=1}^{\infty}$ and their parameters. 
 Define \begin{equation}\label{d-GE}
d_{\text{ge}, j}(t,x,\Omega) := - \left[\partial_t + \Omega\cdot\nabla_x\right] \psi_j + \mathcal L(\psi_j). \end{equation}
Integrate $ | d_{\text{ge},j} |^2$ over $K$, which is equivalent to
\begin{equation}\label{GE1} 
||\partial_t ( \psi - \psi_j) + \Omega\cdot \nabla_x ( \psi - \psi_j) - \mathcal L(\psi) + \mathcal L(\psi_j) ||^2_{L^2(K)}, 
\end{equation}
in which the first two terms are bounded by
\begin{equation}\label{bdd1}
    ||\partial_t ( \psi - \psi_j) + \Omega\cdot \nabla_x ( \psi - \psi_j)||^2_{L^2(K)}
    \lesssim \epsilon_j^2\, T\, 4\pi\, |\mathrm{X}|\, (1 + 16\pi^2), 
\end{equation} 
{where we used the fact that $L^2$ norm is bounded by $L^{\infty}$ norm, due to the boundedness of physical domain $K$. }
The last two terms in \eqref{GE1}, thanks to the operator $\mathcal L$ being bounded, are bounded as follows:
\begin{equation}\label{bdd2}
    ||\mathcal L(\psi) - \mathcal L(\psi_j) ||^2_{{L^2(\mathrm{X}\times\mathbb{S}^2)}}
\lesssim ||\psi - \psi_j ||^2_{{L^2(\mathrm{X}\times\mathbb{S}^2)}} \leq
||\psi - \psi_j ||^2_{{L^{\infty}(\mathrm{X}\times\mathbb{S}^2)}}\,4\pi\, |\mathrm{X}|
< C\, \epsilon_j^2. 
\end{equation} 
Thus the loss term $Loss_{\text{GE}}$ is bounded by $O\left(\epsilon_j^2\right)$. 

For the inflow boundary data, $Loss_{\text{BC}}$ is bounded by 
\begin{equation}\label{bdd3}
    ||\psi_j - \psi|| _{L^2(\gamma_{T}^{-})}^2  \leq T\, 4\pi\, |\partial\mathrm{X}|\, ||\psi_j - \psi||_{L^{\infty}(\gamma_{T}^{-})}^2 \leq T\, 4\pi\, |\partial\mathrm{X}|\, ||\psi_j - \psi||_{L^{\infty}(K)}^2 \leq O\left(\epsilon_j^2\right), 
\end{equation} 
where $\gamma_{T}^{\pm} := [0,T]\times\gamma_{\pm}$. Note that a similar result holds for the specular reflection boundary condition \eqref{rbc}. 

%
For the initial data, we denote the DNN approximation by $\psi_j(0,x,\Omega)$, then 
\begin{equation}\label{bdd4}
    Loss_{\text{IC}} = || \psi_j(0,x,\Omega) - \psi_0(x,\Omega)||_{L^2(\mathrm{X}\times\mathbb{S}^2)}
\leq || \psi_j(0,x,\Omega) - \psi_0(x,\Omega) ||_{L^{\infty}(\mathrm{X}\times\mathbb{S}^2)}\, 4\pi\, |\mathrm{X}| \leq 
O(\epsilon_j^2).
\end{equation} 
Combining all these upper bounds for loss terms in \eqref{bdd1},\eqref{bdd2},\eqref{bdd3} and \eqref{bdd4}, we conclude that
\begin{equation}\label{bdd5}
    Loss_{\text{Total}}(\psi_j) \leq O\left(\frac{1}{j^2}\right).
\end{equation} 
Therefore, $Loss_{\text{Total}}(\psi_j)\to 0$ as $j\to \infty$.

\subsection{Convergence of the DNN solution} 

Now we prove that with the parameters $\{m_{[j]}, w_{[j]}, b_{[j]}\}_{j=1}^{\infty}$ equipped, the neural network in Theorem \ref{Thm1} indeed converges to the analytic solution of the linear transport model, as one desires. 
\begin{theorem}
\label{Thm2}
Let $\{m_{[j]}, w_{[j]}, b_{[j]}\}_{j=1}^{\infty}$ be a sequence defined in Theorem \ref{Thm1}, and $\psi$ is the solution to the linear transport model \eqref{eqn:3Dtransport}. Then, $Loss_{\text{Total}}(\psi_j)\to 0$ implies that
\begin{equation}\label{eqn:thm2}
    ||\psi_j(\cdot, \cdot, \cdot, m_{[j]}, w_{[j]}, b_{[j]}) - \psi ||_{L^{\infty}\left([0,T]; L^2(\mathrm{X}\times\mathbb S^2 )\right)}\to 0. 
\end{equation}  
\end{theorem}
Recall the definition \eqref{d-GE}, since $\psi$ solves the linear transport model \eqref{LTE}, thus
\begin{equation}\label{GE2} \left[\partial_t + \Omega \cdot\nabla_x \right] \left(\psi-\psi_j\right) = d_{\text{ge},j}(t,x,\Omega) + \mathcal L(\psi) - \mathcal L(\psi_j). \end{equation}
Define 
\begin{equation}\label{eqn1}
    d_{\text{ic},j}(x,\Omega) := \psi_0(x,\Omega) - \psi_j(0,x,\Omega), \quad \text{ for } (x,\Omega)\in \mathrm{X}\times\mathbb S^2,
\end{equation} 
in addition to 
\begin{equation}
    d_{\text{bc},j}(t,x,\Omega) := h(t,x,\Omega)- \psi_j(t,x,\Omega)\quad \text{ for }\, (t,x,\Omega)\in\gamma^{-}_{T},
\end{equation} 
for inflow boundary condition \eqref{bc}. 

For simplicity of notations, the $L^2$ norms and the corresponding inner products below stand for in the physical space
$(x,\Omega)\in \mathrm{X}\times\mathbb{S}^2$. Multiplying $(\psi-\psi_j)$ onto \eqref{GE2} and integrating over $\mathrm{X}\times\mathbb{S}^2$, one gets
\begin{align}
\label{GE3}
\begin{split}
&\displaystyle\quad \int_{\mathrm{X}}\int_{\mathbb{S}^2}\partial_t (\psi-\psi_j)^2\, d\Omega dx +
2\int_{\gamma^{+}}(\psi-\psi_j)^2\, \Omega\cdot n_x\, ds - 
2\int_{\gamma^{-}}d_{\text{bc},j}^2\, |\Omega\cdot n_x| \, ds \\[4pt]
&\displaystyle = 2\big\langle \mathcal L(\psi - \psi_j), \psi-\psi_j \big\rangle_{L^2} + 2 \big\langle d_{\text{ge},j}, \psi-\psi_j \big\rangle_{L^2}. 
\end{split}
\end{align}
The first term on the RHS above is bounded by
\begin{equation}  
    2 \big\langle \mathcal L(\psi - \psi_j), \psi-\psi_j \big\rangle_{L^2} \leq \| \mathcal L(\psi - \psi_j)\|_{L^2}^2 + \| \psi- \psi_j \|_{L^2}^2  \lesssim \| \psi - \psi_j \|_{L^2}^2.  
\end{equation} 
The fact that $\int_{\gamma^{+}}(\psi-\psi_j)^2\, \Omega\cdot n_x\, ds \geq 0$ yields
\begin{equation} 
    \frac{d}{dt}||\psi-\psi_j||_{L^2}^2 \leq \underbrace{2\int_{\gamma^{-}}d_{\text{bc},j}^2\, |\Omega\cdot n_x| \, ds + ||d_{\text{ge},j}||_{L^2}^2}_{:=H(t)} + C_1 ||\psi-\psi_j||_{L^2}^2. 
\end{equation} 
Due to $\int_{\gamma^{-}}d_{\text{bc},j}^2\, |\Omega\cdot n_x| \, ds \leq \int_{\gamma^{-}}d_{\text{bc},j}^2\, ds$ and the definitions
\begin{equation} 
    \int_0^t  \int_{\gamma^{-}}\left(d_{\text{bc},j}(\tau,\cdot, \cdot, \cdot)\right)^2\, ds d\tau = Loss_{\text{BC}}, \quad 
    \int_0^t  ||d_{\text{ge},j}(\tau, \cdot, \cdot, \cdot)||_{L^2}^2\, d\tau = Loss_{\text{GE}}, 
\end{equation} 
we have, by Gr\"onwall's inequality, 
\begin{align}
\notag\displaystyle ||\psi(t,\cdot,\cdot)-\psi_j(t,\cdot,\cdot)||_{L^2(\mathrm{X}\times\mathbb{S}^2)}^2 &\leq e^{C_1 t } Loss_{\text{IC}} + e^{C_1 t} \int_0^t H(\tau) d\tau \\[4pt]
&\notag\displaystyle \lesssim e^{C_1 t}\left(Loss_{\text{IC}} + 2 Loss_{\text{BC}} + Loss_{\text{GE}}\right) \\[4pt]
&\displaystyle \leq 2 e^{C_1 t}\, Loss_{\text{Total}}. 
\end{align}

It's already known in \eqref{bdd5} from Theorem \ref{Thm1} that $Loss_{\text{Total}}(\psi_j) \leq O(\frac{1}{j^2})$, thus
\begin{equation}\label{R1}
    ||\psi(t,\cdot,\cdot)-\psi_j(t,\cdot,\cdot)||_{L^2(\mathrm{X}\times\mathbb{S}^2)}^2 \leq 
    C^{\prime}\, e^{C_1 T}\,\frac{1}{j^2}, \quad \text{ for } t \in [0,T].
\end{equation}
After taking $L^{\infty}$ norm in $t\in [0,T]$, we conclude that 
\begin{equation} 
    ||\psi-\psi_j||_{L^{\infty}([0,T]; L^2(\mathrm{X}\times\mathbb{S}^2))} \to 0, \quad \text{as   }\, j\to \infty\,. 
\end{equation} 

\begin{remark}
We remark that both Theorem \ref{Thm1} and Theorem \ref{Thm2} hold true for the case of specular reflection boundary condition \eqref{rbc} with corresponding loss term defined in \eqref{Loss_Spec}, and the proof is similar thus omitted here. 
\end{remark}

\section{Numerical Experiment}
\label{sec:Num}

\subsection{Grid points}

To approximate the distribution function by the DNN approach, we make the data of grid points for each variable. 
The grid points of $t$, $x$ and $\mu$ for the training are chosen uniformly as follows: 
\begin{equation}
    \{ (t_i, x_j, \mu_k) \}_{i,j,k} \in [0,T]\times\Omega_x\times\Omega_{\mu}, \quad\text{  with fixed  } \,\Delta t,\, \Delta x,\,\Delta\mu.  
\end{equation} 


The sampling points for Loss$_\text{BC}$ and Loss$_\text{IC}$ are chosen as random points.
\subsection{Numerical Tests}
{In this section, we present simulation results for several benchmark problems of the model equation.  }
By assuming slab geometry, equation \eqref{eqn:3Dtransport} reduces to a one-dimensional problem:
\begin{equation}\label{eqn:slab_transport}
    \frac{1}{c}\frac{\partial \psi}{\partial t} + \mu \partial_x \psi + \sigma_t \psi = \sigma_s \mathcal{S} \psi + \mathcal{Q},
\end{equation}
where $x \in (x_L,x_R)$ is a scalar coordinate perpendicular to the slab, and $\mu \in[-1,1]$ is the cosine of the angle between $\Omega$ and the positive x-axis. 
{Then the scattering operator is simplified as
\begin{equation}
(\mathcal{S} \psi)(x)  = \frac{1}{2}\int_{-1}^1 \, g(\mu) \psi(x,\mu) \, d\mu.
\end{equation}
For the one-dimensional problem with isotropic scattering, the kernel function $g \equiv 1$ is a constant.}  { We shall compute the following integral in the DNN by Gaussian quadrature as:
\begin{equation}
{(\mathcal{S} \psi)(x)   = \frac{1}{2}\int_{-1}^1 \, \psi(x,\mu) \,  d\mu \simeq } \frac{1}{2}\sum_{j=1}^{N_{\mu}}\,\xi_j\,\psi(t,x,\mu_j),
\end{equation}
which discrete along the $\mu$ direction and $\xi_j$ denotes the corresponding weights. In the following numerical experiments, we shall also use the notations: angular average $\rho = \frac{1}{2}\int_{-1}^1\,\psi \, d\mu$ and the corresponding discrete approximation $\rho_h = \frac{1}{2}\sum_{j=1}^{N_\mu}\xi_j\,\psi(t,x,\mu_j).$
}

In the following approximation, we shall choose 8 layers including 1 input layer and 1 output layer, and each hidden layer with 64 neurons. Besides, we choose $100$ points uniformly along $x$ and $64$ Gaussian quadrature along $\mu$ directions; $50$ points uniformly along $t$ direction; {we choose 200 random points for approximating the boundary condition and initial condition.}

\subsubsection{Manufactured solution}\label{Sect:Num_Test0}

\begin{figure}[H]
    \centering
    \begin{tabular}{cc}
    \includegraphics[width=.45\textwidth]{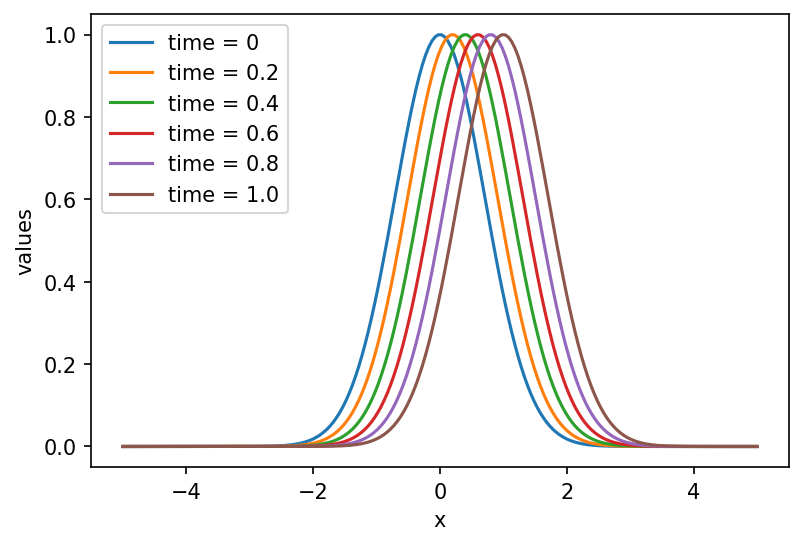}
    &\includegraphics[width=.45\textwidth]{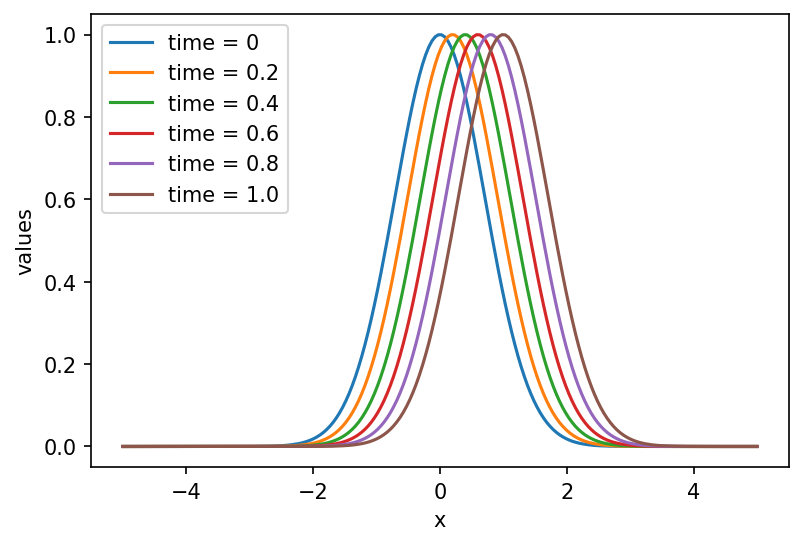}\\
    (a) & (b)
    \end{tabular}
    \caption{Example \ref{Sect:Num_Test0} with $\sigma_s = 0.0$ and $\sigma_a = 1.0$: (a). exact solution $\rho$; (b). DNN approximation $\rho_h$.}
    \label{fig:Test0-S0A1}
\end{figure}
In this example, we set $x\in(-5,5)$, $\mu\in [-1,1]$, and $t\in(0,1)$ and the exact solution is taken as
\begin{equation}
    \psi = \exp(-(x-t)^2).
\end{equation}
Then the source term are chosen to fit the exact solution and parameter values in $\sigma_s$ and $\sigma_a$.

\vskip.1in

\begin{figure}[H]
\centering
\includegraphics[width=.7\textwidth]{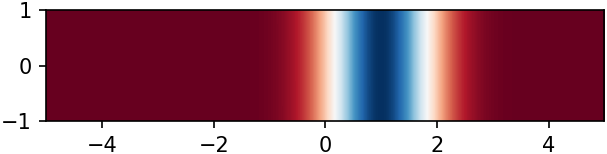}
\caption{Example \ref{Sect:Num_Test0} with $\sigma_s = 0.0,\sigma_a = 1.0$: DNN approximation to $\psi$ at time = 1.0.}
\label{Fig:Test0-S0A1-2D-t1}
\end{figure}

First, as we choose $\sigma_s = 0.0$ and {$\sigma_a = 1.0$} for a purely absorbing case, the exact solution and DNN approximation for $\rho$ are plotted in Figure \ref{fig:Test0-S0A1} (a)-(b). It is observed that the simulations induced by DNN match the exact solution at different time period. The DNN solution at time$= 1.0$ on the $x-\mu$ plane is plotted in Figure \ref{Fig:Test0-S0A1-2D-t1}. As we can observe the uniform distribution along the $\mu$ direction.

\begin{figure}[H]
    \centering
    \begin{tabular}{cc}
    \includegraphics[width=.45\textwidth]{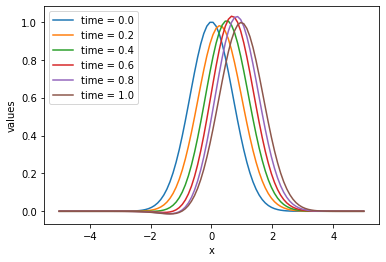}
    &\includegraphics[width=.45\textwidth]{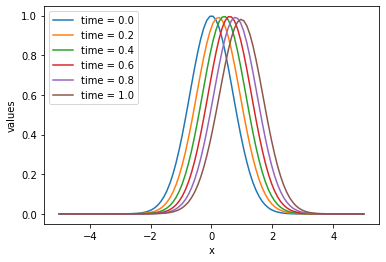}\\
    (a) & (b)
    \end{tabular}
    \caption{Example \ref{Sect:Num_Test0}: Plot of DNN solutions for $\rho_h$: (a) $\sigma_s = 1$ and $\sigma_a = 0$; (b) $\sigma_s = 1.0$ and $\sigma_a = 1.0$.}
    \label{fig:Test0-S1A0A10}
\end{figure}

\begin{figure}[!ht]
\centering
\includegraphics[width=.7\textwidth]{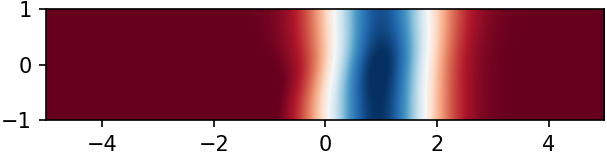}
\caption{Example \ref{Sect:Num_Test0} with $\sigma_s = 1.0,\sigma_a = 0.0$: DNN approximation to $\psi$ at time = 1.0.}
\label{Fig:Test0-S1A0-2D-t1}
\end{figure}

\begin{figure}[!ht]
\centering
\includegraphics[width=.7\textwidth]{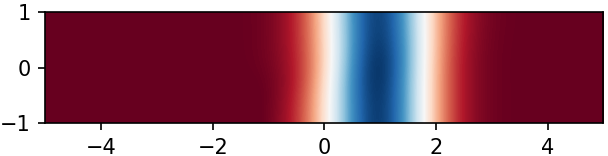}
\caption{Example \ref{Sect:Num_Test0} with $\sigma_s = 1.0,\sigma_a = 1.0$: DNN approximation to $\psi$ at time = 1.0.}
\label{Fig:Test0-S1A1-2D-t1}
\end{figure}

Next, we choose $\sigma_s = 1.0$ and choose different values in {$\sigma_a$}. The DNN solutions {for $\rho_h$} corresponding to { $\sigma_a = 0.0$ (purely scattering case) and $\sigma_a = 1.0$} are plotted in Figure \ref{fig:Test0-S1A0A10}. We observe that: in the purely scattering case, the DNN solutions show larger error as comparing to the exact solution. Then the 2-dimensional plots for DNN solutions at time = 1.0 on the $x-\mu$ plane are plotted in Figure \ref{Fig:Test0-S1A0-2D-t1}-\ref{Fig:Test0-S1A1-2D-t1}.

Lastly, the relative errors corresponding to different parameter settings for various time periods are reported in Table \ref{tab:Test0}. It shows that the purely absorbing case shows the least error.  {This is not surprising, since the collisions (through the scattering operator) introduce more micro features and thus makes the particle dynamic more complicated. }

\begin{table}[H]
    \centering
    \caption{Example \ref{Sect:Num_Test0}: Relative Errors in DNN approximations to $\psi$.}
    \label{tab:Test0}
    \begin{tabular}{c|c|c|c}\hline\hline
    time& $\sigma_s = 0.0,\sigma_a = 1.0$ & $\sigma_s = 1.0,\sigma_a = 0.0$ & $\sigma_s = 1.0,\sigma_a = 1.0$ \\ \hline
    0.0 &2.3329e-4 &2.7779e-4 &2.8335e-4\\
    0.2 &1.8372e-4 &4.4767e-2 &2.1264e-2\\
    0.4 &1.5604e-4 &4.1943e-2 &2.7555e-2\\
    0.6 &1.3525e-4 &4.1276e-2 &2.6886e-2\\
    0.8 &1.1590e-4 &4.7595e-2 &2.8699e-2\\
    1.0 &1.1543e-4 &4.5483e-2 &2.1431e-2\\ \hline\hline
    \end{tabular}
\end{table}

\subsubsection{The plane source}\label{sect:Num_Test1}
Initially some particles are emitted from a planar source to infinite medium.  Again thanks to the symmetry of the problem, the model is allowed to be reduced to a one-dimensional spatial setting, that is, using variable $x$ to measure the signed normal distance from the location to the source plane. 

\begin{figure}[H]
    \centering
    \begin{tabular}{cc}
    \includegraphics[width=0.45\textwidth]{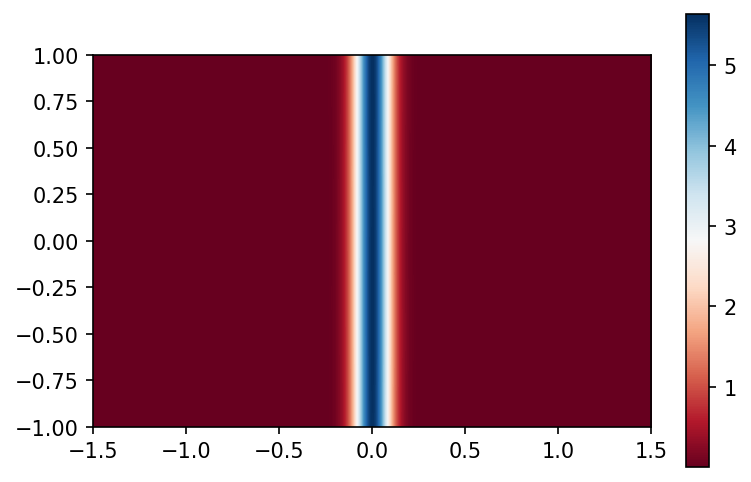}
&\includegraphics[width=0.45\textwidth]{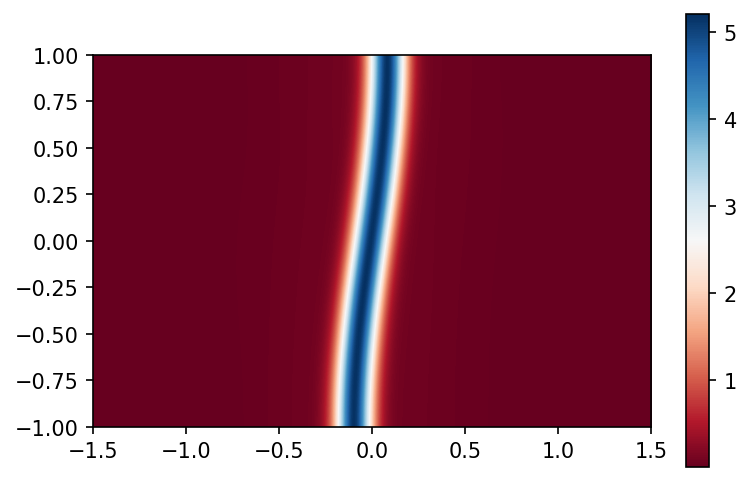}\\
(a) time = 0.0 & (b) time = 0.1\\
\includegraphics[width=0.45\textwidth]{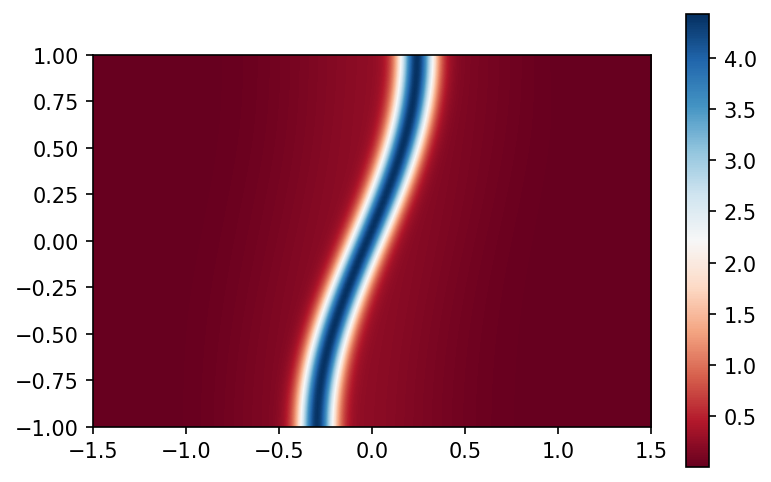}
&\includegraphics[width=0.45\textwidth]{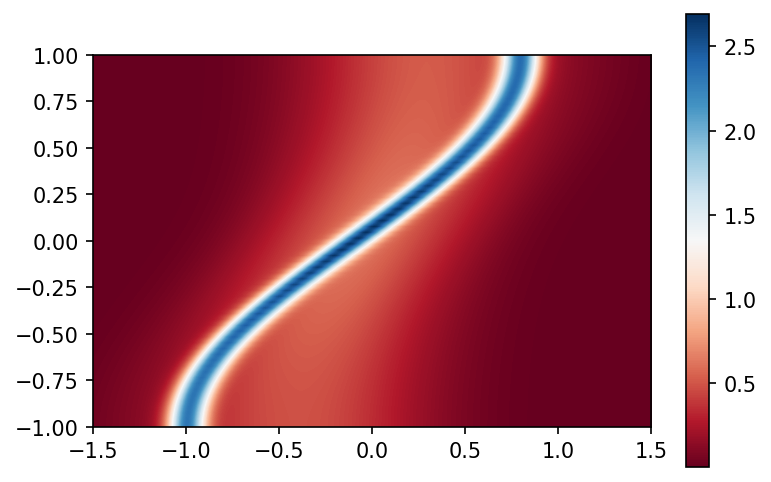}
\\
(c) time = 0.3 & (d) time = 1.0
\end{tabular}
    \caption{Example \ref{sect:Num_Test1}: Plot of DNN approximation to $\psi$ for $\sigma_s = 1.0$ and $\sigma_a = 0.0$ at different time. }
    \label{fig:Test1-S1A0}
\end{figure}

\begin{figure}
    \centering
    \begin{tabular}{cc}
    \includegraphics[width=0.45\textwidth]{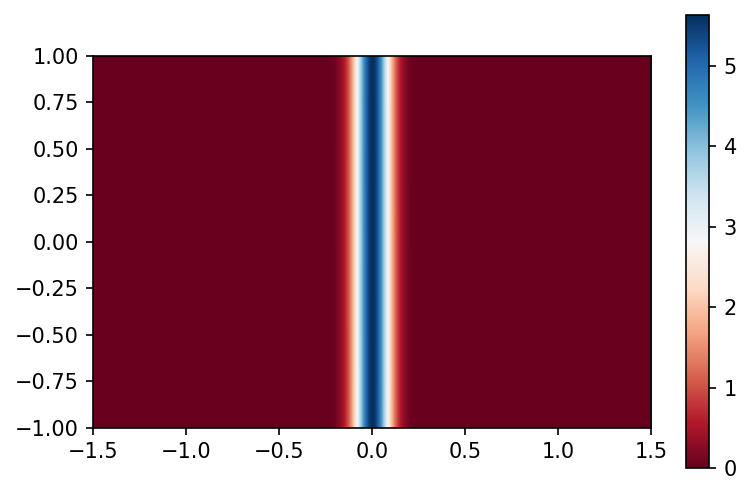}
&\includegraphics[width=0.45\textwidth]{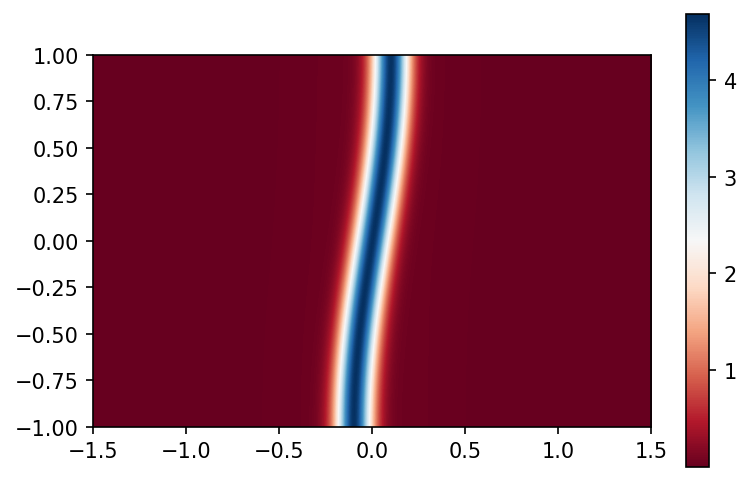}\\
(a) time = 0.0 & (b) time = 0.1\\
\includegraphics[width=0.45\textwidth]{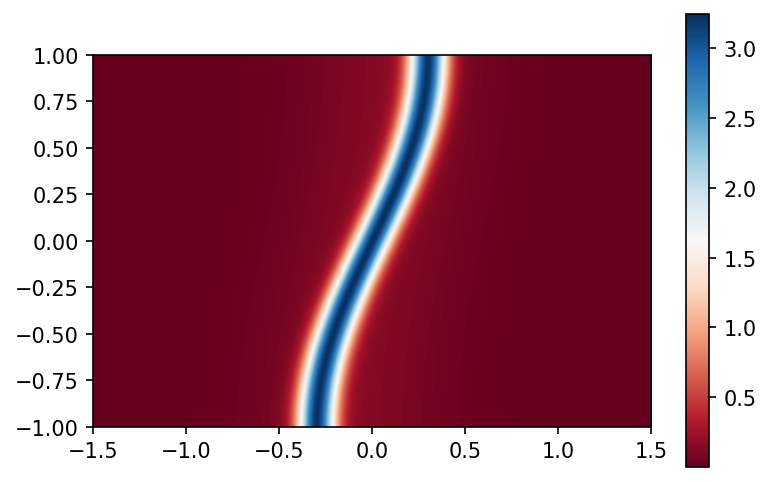}
&\includegraphics[width=0.45\textwidth]{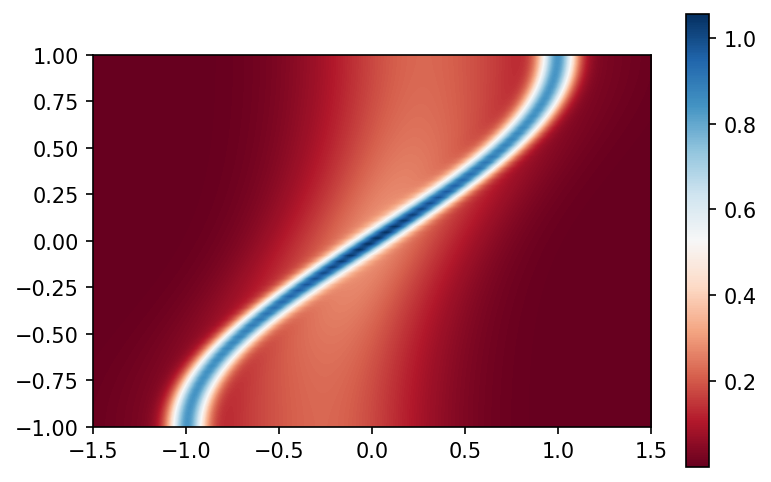}
\\
(c) time = 0.3 & (d) time = 1.0
\end{tabular}
    \caption{Example \ref{sect:Num_Test1}: Plot of DNN approximation to $\psi$ for $\sigma_s = 1.0$ and $\sigma_a = 1.0$ at different time. }
    \label{fig:Test1-S1A1}
\end{figure}

\begin{figure}
\begin{tabular}{ccc}
\includegraphics[width=0.45\textwidth]{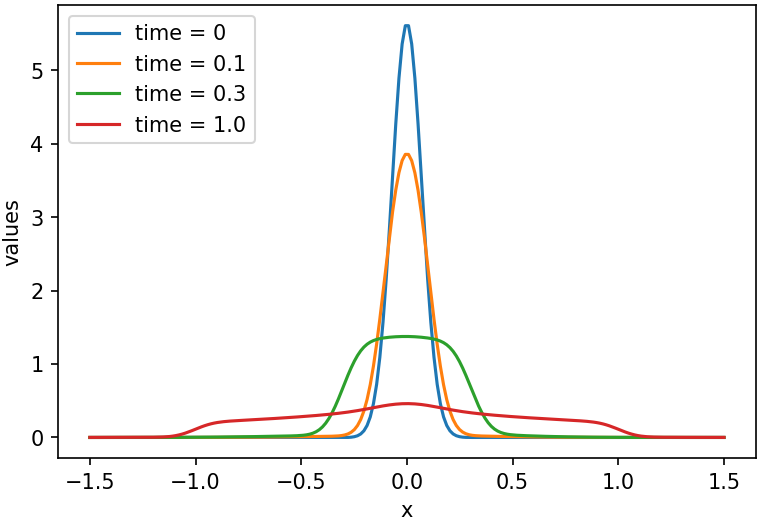}
&\includegraphics[width=0.45\textwidth]{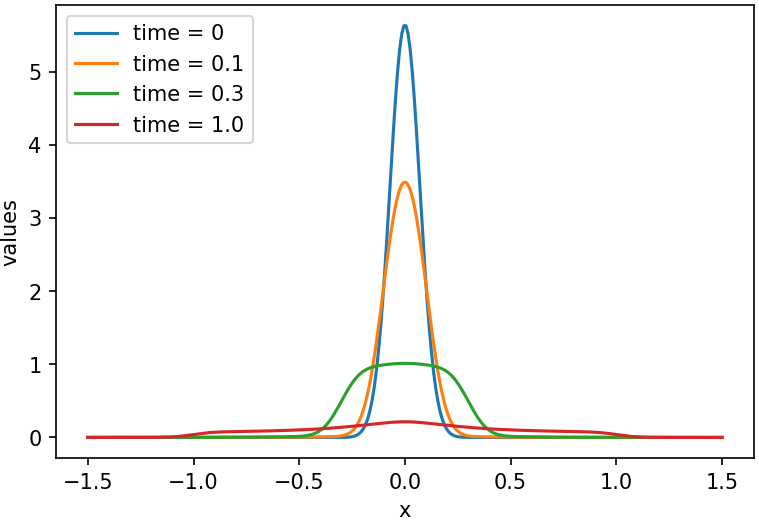}\\
(a) & (b) 
\end{tabular}
\caption{Example \ref{sect:Num_Test1}: Plots of DNN solutions $\rho_h$ for $k = 100$ and (a) $\sigma_s = 1.0$, $\sigma_a = 0.0$; (b)$\sigma_s = 1.0$, $\sigma_a = 1.0$.}\label{Fig:Test1-S1A0}
\end{figure}

{
\vskip.1in
In this test, we assume there is no external source (i.e. $\mathcal{Q} = 0$) and particle speed $c=1$, and choose domain $(x,\mu) \in (-1.5,1.5)\times[-1,1]$ and the initial condition as below with $k = 100$:
\begin{equation}\label{eqn:ic_gaussian}
\psi|_{t = 0} = \exp(-kx^2)\frac{\sqrt{k}}{\text{erf}(\sqrt{k})\sqrt{\pi}}.
\end{equation}

Again we test on different values of cross-sections:
\begin{enumerate}
    \item $\sigma_s = 1.0$ and $\sigma_a = 0.0$: the approximation calculated by DNN is shown in Figure \ref{fig:Test1-S1A0} for time $t=0,0.1,0.3,1.0$. Then after we take average along $\mu$ direction, the 1-dimensional plot for $\rho$ is shown in Figure \ref{Fig:Test1-S1A0} (a). As one can observe from this figure that since the absorbing coefficient $\sigma_a$ is zero, the system shows the purely scattering.
    \item  $\sigma_s = 1.0$ and $\sigma_a = 1.0$: the approximation calculated by DNN is shown in Figure \ref{fig:Test1-S1A1}. Then after we take average along $\mu$ direction, the 1-dimensional plot is shown in Figure \ref{Fig:Test1-S1A0} (b). Due to the appearing of absorbing, although the solution at time = 1.0 shows the similar shape as that in Figure \ref{Fig:Test1-S1A0} (a), the magnitude in Figure \ref{Fig:Test1-S1A0} (b) is significantly reduced. {Note that the level bars of color-map for Figure \ref{fig:Test1-S1A0} and Figure \ref{fig:Test1-S1A1} are not the same.}
\end{enumerate}


}


\begin{figure}
    \centering
    \includegraphics[width=0.45\textwidth]{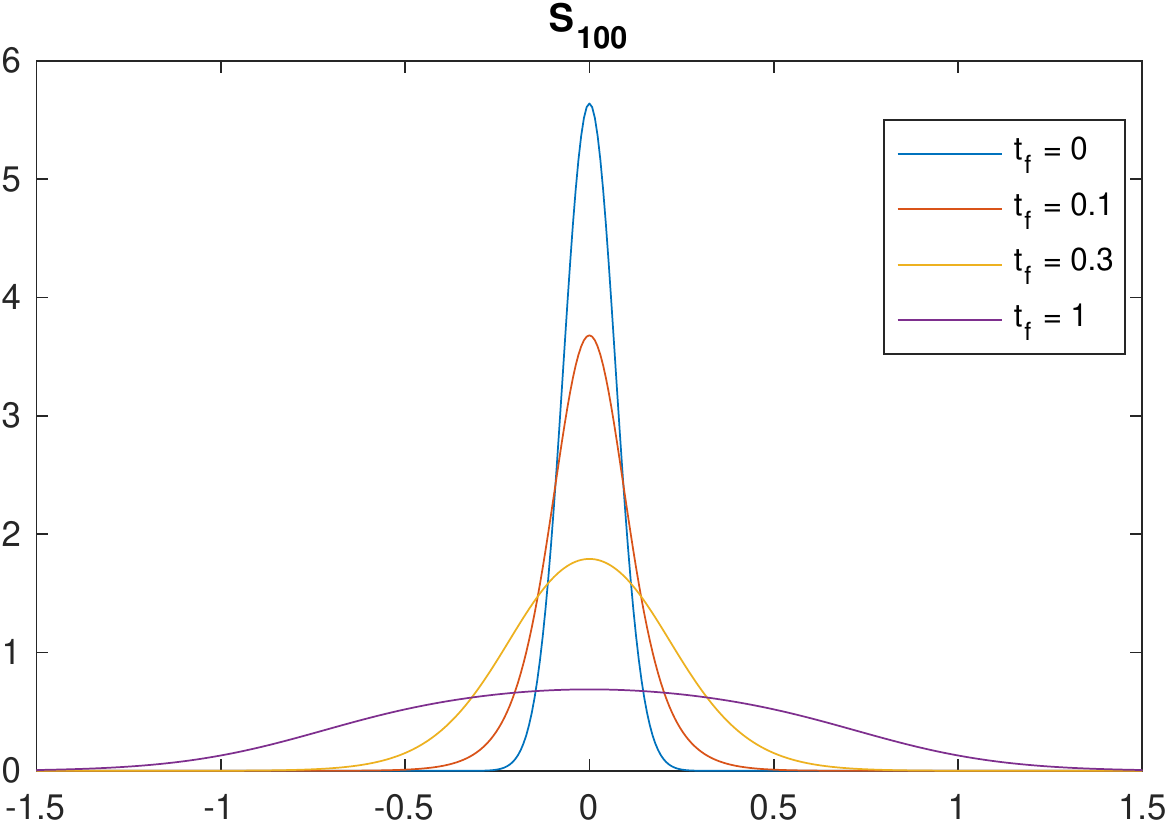}
    \includegraphics[width=0.45\textwidth]{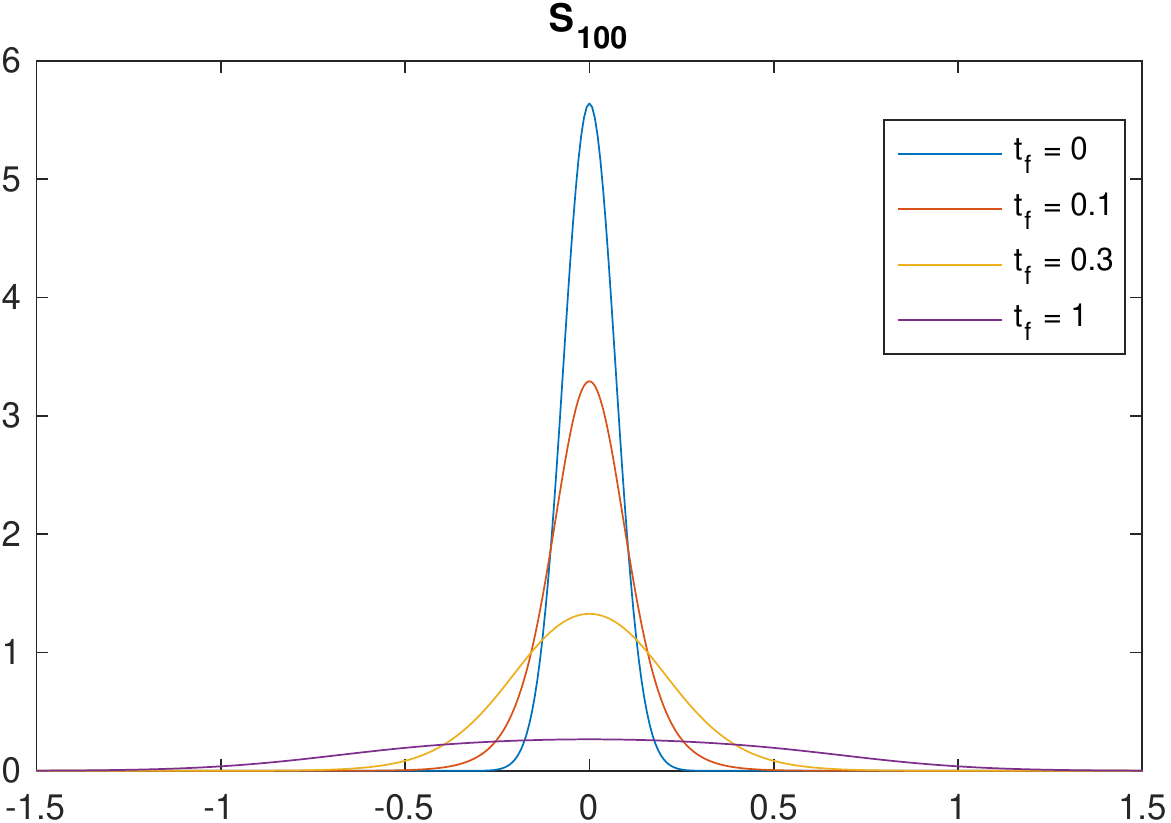}
    \caption{Example \ref{sect:Num_Test1}: Plots of angular average of discrete ordinate $S_{100}$ solutions
    for 
    $k=100$: 
     (a) $\sigma_s = 1$, $\sigma_a = 0$;  (b) $\sigma_s = 1$, $\sigma_a = 1$.
    }
    \label{fig:Test1_ref_SN}
\end{figure}

{The angular average of numerical simulations with discrete ordinate methods \cite{hauck2013collision} $S_N$ with $N=100$ are presented in Figure \ref{fig:Test1_ref_SN}.  Comparing Figure \ref{fig:Test1_ref_SN} with Figure \ref{Fig:Test1-S1A0}, it seems that the DNN solutions have a more flat peak at $t=0.3$ and spreads out faster than $S_N$ solutions. We do not yet understand the origin of this behavior and {will }investigate more in the future. }
\subsubsection{Two beam problem}\label{Sect:Num-3}
\begin{figure}
\centering
\includegraphics[width=0.33\textwidth]{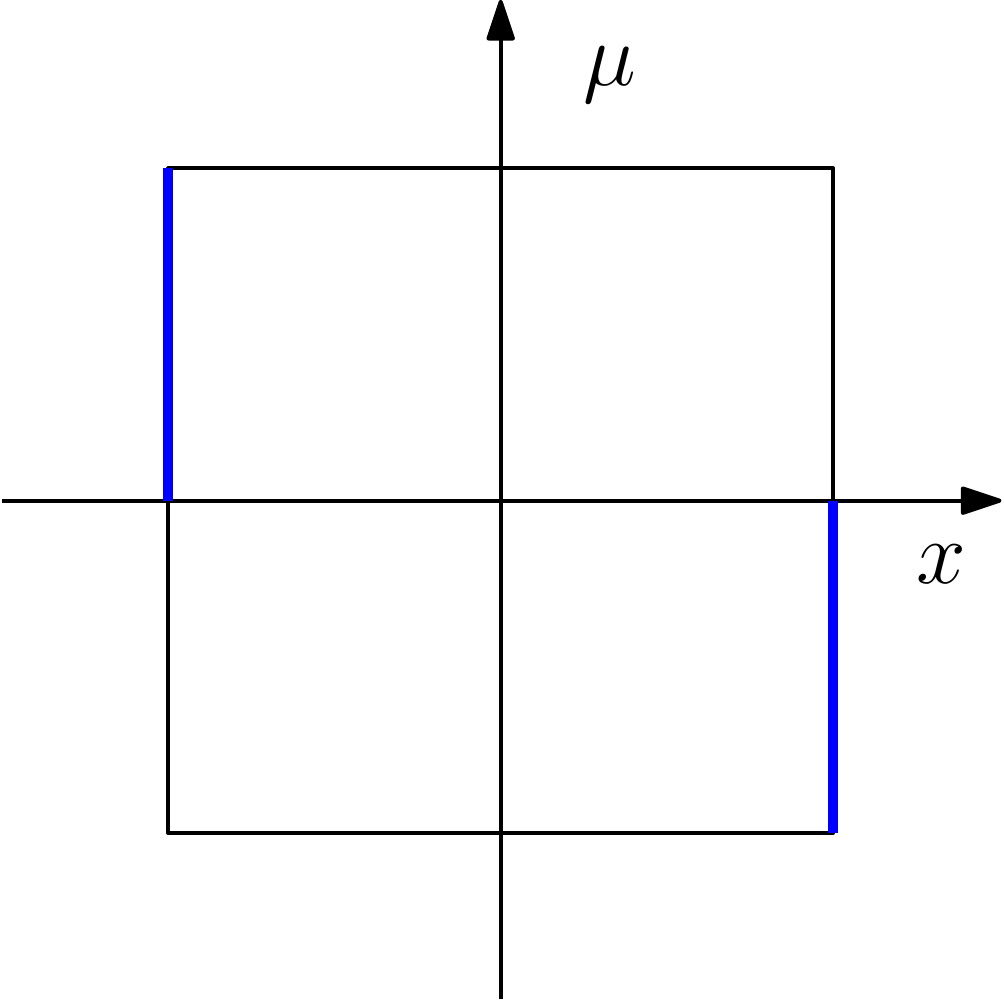}
\caption{Example \ref{Sect:Num-3}: The illustration of the boundary condition.}\label{Fig:Test2-bc}
\end{figure}

In this problem, particles enter the material slab from both boundaries. The spatial domain is $x\in [-1,1]$, and we assume isotropic scattering and $c=1$, $\mathcal{Q} = 0$
. Initially, it is a void, that is, $\psi(t=0,x,\mu) = 0$.
Boundary conditions are isotropic incoming fluxes at both boundaries:
\begin{equation}
    \psi(t,-1,\mu>0) = \psi(t,1,\mu<0) = 10.
\end{equation}
The boundary condition is also illustrated in Figure \ref{Fig:Test2-bc}.

\begin{figure}[H]
\begin{tabular}{cc}
\includegraphics[width=0.45\textwidth,height = .35\textwidth]{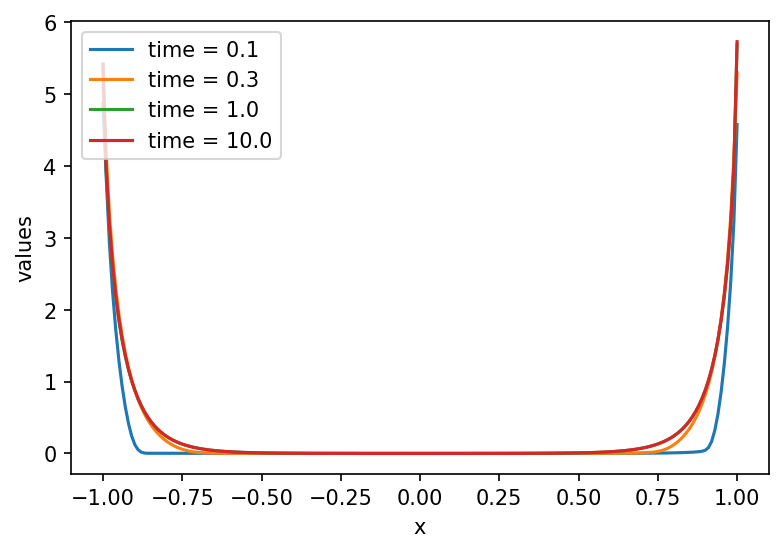}
&\includegraphics[width=0.45\textwidth]{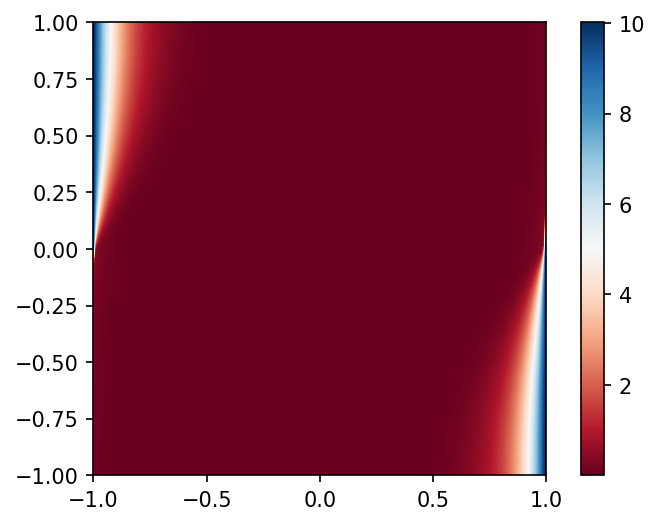}\\
(a) & (b)
\end{tabular}
\caption{Example \ref{Sect:Num-3}: Case (1) $\sigma_s = 1.0$, $\sigma_a = 9.0$: (a). solution of $\rho_h$ at different time; (b). 2-dimensional plot of DNN approximation to $\psi$ on the $x-\mu$ plane at $t=10$.}\label{Fig:Test2-S1A9}
\end{figure}

\begin{figure}[H]
\begin{tabular}{cc}
\includegraphics[width=0.45\textwidth,height = .35\textwidth]{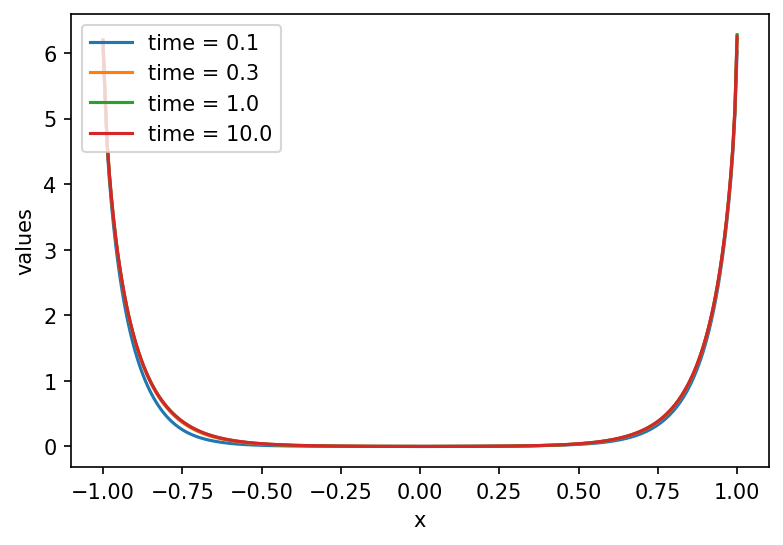}
&\includegraphics[width=0.45\textwidth]{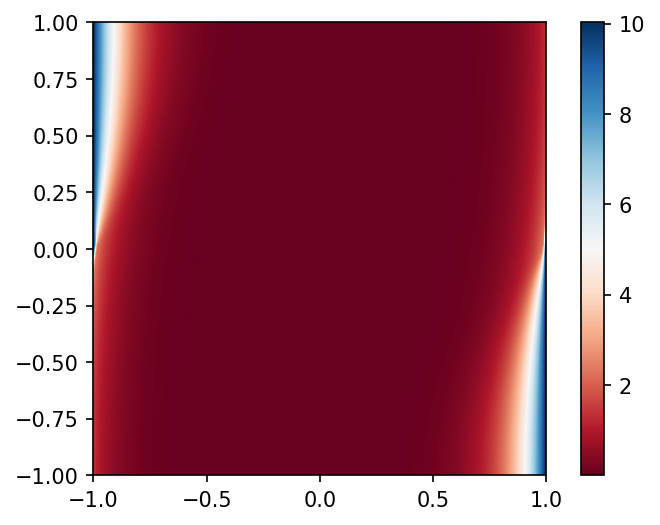}\\
(a) & (b)
\end{tabular}
\caption{Example \ref{Sect:Num-3}: Case (2) $\sigma_s = 5.0$, $\sigma_a = 5.0$: (a). solution of $\rho_h$ at different time; (b). 2-dimensional plot of DNN approximation to $\psi$ on the $x-\mu$ plane at $t=10$.}\label{Fig:Test2-S5A5}
\end{figure}

\begin{figure}[H]
\begin{tabular}{cc}
\includegraphics[width=0.45\textwidth,height = .35\textwidth]{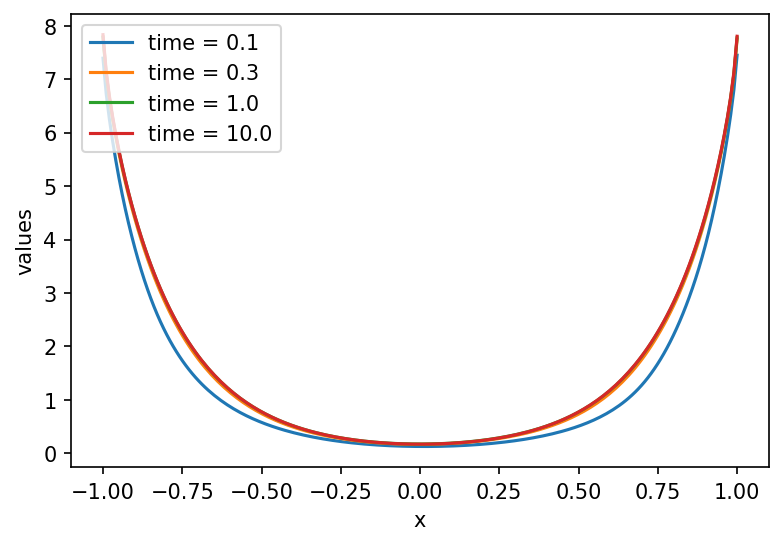}
&\includegraphics[width=0.45\textwidth]{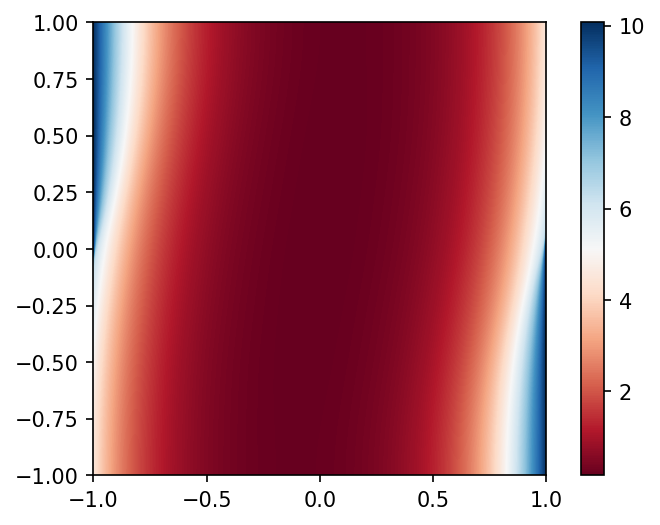}\\
(a) & (b)
\end{tabular}
\caption{Example \ref{Sect:Num-3}: Case (3) $\sigma_s = 9.0$, $\sigma_a = 1.0$: (a). solution of $\rho_h$ at different time; (b). 2-dimensional plot of DNN approximation to $\psi$ on the $x-\mu$ plane at $t=10$.}\label{Fig:Test2-S9A1}
\end{figure}

\begin{figure}[H]
\resizebox{\textwidth}{!}{
\begin{tabular}{ccc}
\includegraphics[width=0.33\textwidth]{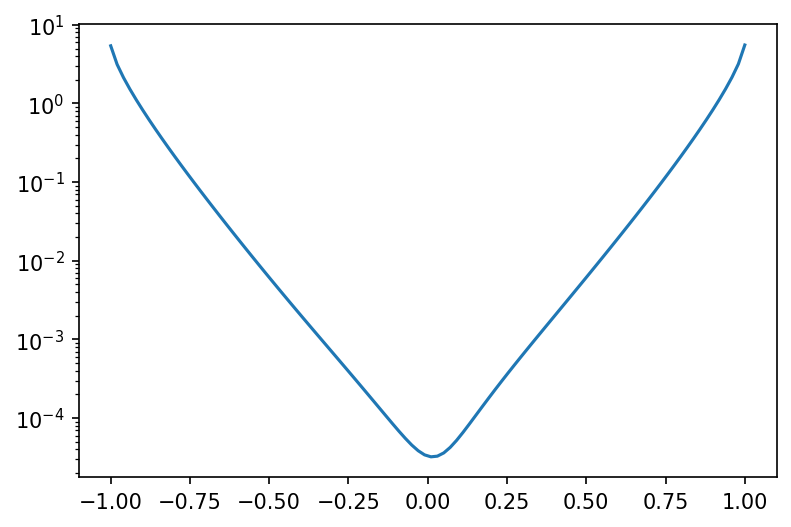}
&\includegraphics[width=0.33\textwidth]{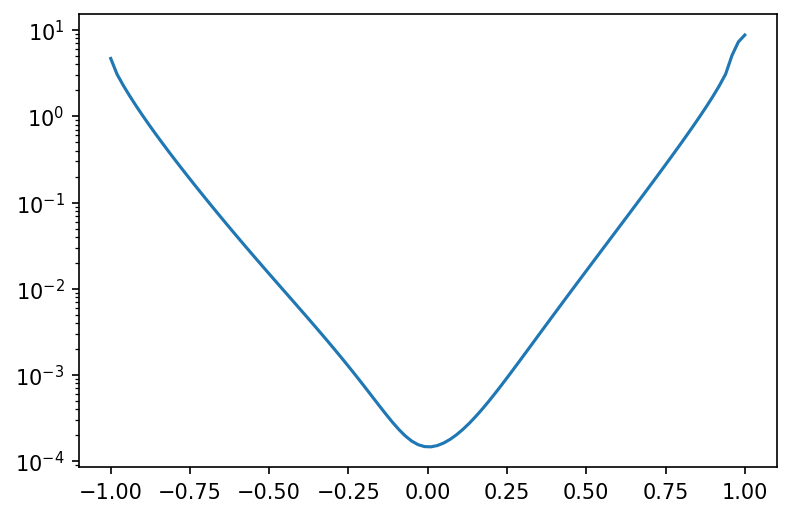}
&\includegraphics[width=0.33\textwidth]{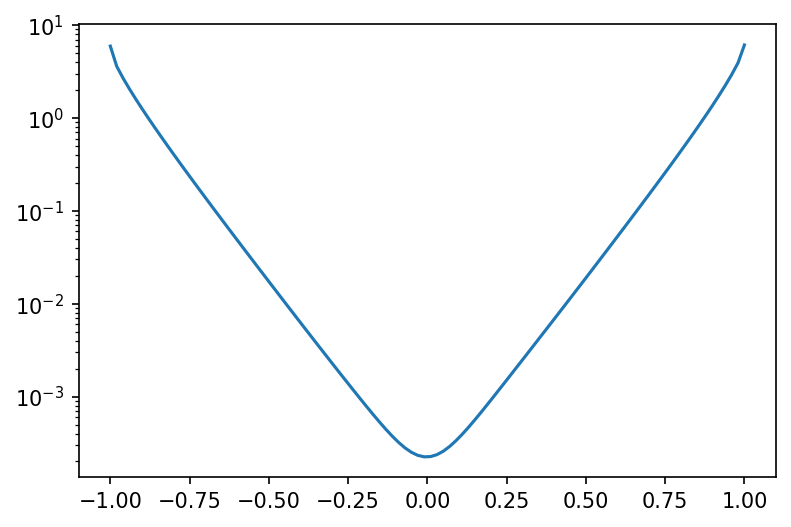}\\
(1) $\sigma_s = 1.0$, $\sigma_a = 9.0$ &
(2) $\sigma_s = 5.0$, $\sigma_a = 5.0$ & 
(3) $\sigma_s = 9.0$, $\sigma_a = 1.0$
\end{tabular}
}
\caption{Example \ref{Sect:Num-3}: Plots of DNN solutions of $\rho_h$ at time = 10.0.}
\label{Fig:Test2-t10}
\end{figure}

We consider the following cases with different values of 
scattering ratio $r = \frac{\sigma_s }{\sigma_t}$ with $\sigma_t = 10$ and time up to $t=10$:
\begin{enumerate}
    \item[(1)] absorption dominated case: $r = 0.1$;
    \item[(2)]  balanced case: $r = 0.5$;
    \item[(3)]  scattering dominated case: $r = 0.9$.
\end{enumerate}
For the testing cases, the DNN solutions are plotted in Figure \ref{Fig:Test2-S1A9}-\ref{Fig:Test2-S9A1} and the semi-log plots for time $t=10$ { are shown} 
in Figure \ref{Fig:Test2-t10}.
The time snapshots of angular average show the  process how particles enter through boundaries and get to reach a steady state.  
Comparing the results with different scattering ratios, one observes that the higher the scattering ratio is (i.e. the weaker the absorption is), 
the higher the center valley of angular average reaches and the more total mass is. Such observations meet the physical expectations. 


\section{Conclusions and future work}
\label{sec:Con}
This paper studies the DNN approach in solving linear transport equations. Theoretical analysis is conducted, and we show two main results: 
1) the loss function will go to zero; 
2) the neural network solution will converge point-wisely to the analytic solution. 
Effectiveness and efficiency of the proposed method are validated in the numerical simulations. 

{It is noticed that in Figure \ref{fig:Test0-S1A0A10} (a), around the left shift from void to the peak, the value of density function is negative, although the magnitude is very small, which is also commonly observed in traditional numerical methods using polynomials as approximations for near-void area.  In the future, we are interested in enforcing {positivity of the solution} in the DNN approach to improve or fix this issue.} 
For future work, we  shall extend this DNN approach to study high-dimensional kinetic problems with random parameters and multiple scales. 
Another important and interesting topic in the rarefied gas dynamics is to study the asymptotic behaviors of particles, such as convergence towards the equilibrium state of the distribution and behaviors of macroscopic quantities. We will explore answers to these questions when using the deep learning methods.

\bibliographystyle{siam}
\bibliography{DNN_Ref.bib}

\end{document}